\documentclass[11pt,reqno]{amsart}
\usepackage{amssymb}
\usepackage{amsmath}
\usepackage{amscd}
\begin{document}
\setlength{\oddsidemargin}{0cm} \setlength{\evensidemargin}{0cm}

\theoremstyle{plain}
\newtheorem{theorem}{Theorem}[section]
\newtheorem{proposition}[theorem]{Proposition}
\newtheorem{lemma}[theorem]{Lemma}
\newtheorem{corollary}[theorem]{Corollary}
\newtheorem{conj}[theorem]{Conjecture}

\theoremstyle{definition}
\newtheorem{definition}[theorem]{Definition}
\newtheorem{exam}[theorem]{Example}
\newtheorem{remark}[theorem]{Remark}

\numberwithin{equation}{section}

\title[Invariant Einstein metrics on three-locally-symmetric spaces]
{Invariant Einstein metrics on three-locally-symmetric spaces}

\author{Zhiqi Chen}
\address{School of Mathematical Sciences and LPMC \\ Nankai University \\ Tianjin 300071, P.R. China} \email{chenzhiqi@nankai.edu.cn}

\author{Yifang Kang}
\address{Institute of Mathematics and Physics, Central South University of Forestry and Technology, Changsha Hunan, 410004, P. R. China} \email{kangyf@tju.edu.cn}

\author{Ke Liang}
\address{School of Mathematical Sciences and LPMC \\ Nankai University \\ Tianjin 300071, P.R. China}\email{liangke@nankai.edu.cn}



\begin{abstract}
In this paper, we classify three-locally-symmetric spaces for a connected, compact and simple Lie group. Furthermore, we give the classification of invariant Einstein metrics on these spaces.
\end{abstract}

\maketitle


\setcounter{section}{0}
\section{Introduction}
A Riemannian manifold $(M,\langle\cdot,\cdot\rangle)$ is called Einstein
if the Ricci tensor $\mathrm{Ric}$ of the metric $\langle\cdot,\cdot\rangle$ satisfies $\mathrm{Ric}=c\langle\cdot,\cdot\rangle$ for some
constant $c$. The above Einstein equation reduces to a system of nonlinear second-order partial differential
equations. But it is difficult to get general existence results. Under the assumption that $M$ is a homogeneous Riemannian
manifold, the Einstein equation reduces to a more manageable system of (nonlinear) polynomial equations, which in some cases can be solved explicitly.
There are a lot of progress in the study on invariant Einstein metrics of homogeneous manifolds, such as the articles \cite{Ar1,AK,AC1,Bo1,Bo2,Bo3,DZ1,GLP,He1,Ke1,Ki1,Ja1,LNF1,Mo1,Mu1,Ni00,Po1,Ro1,Sa1,Wa1,WZ1},
and the survey article \cite{NR1} and so on.

Consider a homogeneous compact space $G/H$ with a semisimple connected Lie group $G$ and a connected Lie subgroup $H$. Denote by $\mathfrak g, \mathfrak h$ the Lie algebras of $G,H$ respectively. Assume that $\mathfrak p$ is the orthogonal complement of $\mathfrak h$ in $\mathfrak g$ with respect to $B$, where $B$ is the Killing from of $\mathfrak g$. Every $G$-invariant metric on $G/H$ generates an $\rm{ad}{\mathfrak h}$-invariant inner product on $\mathfrak p$ and vice versa \cite{Be1}. This makes it possible to identify invariant Riemannian metrics on $G/H$ with $\rm{ad}{\mathfrak h}$-invariant inner product on $\mathfrak p$. Note that the metric generated by the inner product $-B|_{\mathfrak p}$ is called standard. Furthermore, if $G$ acts almost effectively on the homogeneous space $G/H$, and $\mathfrak p$ is the direct sum of three $\rm{ad}{\mathfrak h}$-invariant irreducible modules pairwise orthogonal with respect to $B$, i.e.,
$$\mathfrak p=\mathfrak p_1\oplus \mathfrak p_2\oplus\mathfrak p_3,$$
with $[\mathfrak p_i,\mathfrak p_i]\subset \mathfrak h$ for any $i\in \{1,2,3\}$, then $G/H$ is called a three-locally-symmetric space.

The notation of a three-locally-symmetric space is introduced by Nikonorov in \cite{Ni00}. There has been a lot of study on invariant Einstein metrics for certain three-locally-symmetric spaces. For example, invariant Einstein metrics on the flag manifold $SU(3)/T_{max}$ are given in \cite{DN1}, on $$Sp(3)/Sp(1)\times Sp(1)\times Sp(1) \text{ and } F_4/Spin(8)$$ are obtained in \cite{Ro1}, on the K$\rm{\ddot{a}}$hler C-spaces $$SU(n_1+n_2+n_3)/S(U(n_1)\times U(n_2)\times U(n_3)), SO(2n)/U(1)\times U(n-1), E_6/U(1)\times U(1)\times Spin(8)$$ are classified in \cite{Ki1}, another approach to $$SU(n_1+n_2+n_3)/S(U(n_1)\times U(n_2)\times U(n_3))$$ is given in \cite{Ar1}. The existence is proved in \cite{Ni00} of at least one invariant Einstein metric for every three-locally-symmetric space. Furthermore in \cite{LNF1}, the classification of invariant Einstein metrics is given for $$Sp(l+m+n)/Sp(l)\times Sp(m)\times Sp(m)\text{ and } SO(l+m+n)/SO(l)\times SO(m)\times SO(m).$$
But the classification of three-locally-symmetric spaces is still incomplete, which leads to the incomplete classification of invariant Einstein metrics. In this paper, we complete the classification of three-locally-symmetric spaces for $G$ simple, and then classify invariant Einstein metrics.

The paper is organized as follows. In Section 2, we give the correspondence between the classification of three-locally-symmetric spaces $G/H$ and that of certain involution pairs of $G$. In Section 3, the classification of  three-locally-symmetric spaces $G/H$ is given for $G$ simple based on the theory on involutions of compact Lie groups. We list them in Table 1 in Theorem~\ref{tlss}. Furthermore we prove the isotropy summands are pairwise nonisomorphic for the new three-locally-symmetric spaces. It makes the method of \cite{LNF1} valid to give the classification of invariant Einstein metrics on these three-locally-symmetric spaces in Section 4.

\begin{remark}
The minute that we uploaded this paper on www.arXiv.org, we received an Email from Prof. Yurii Nikonorov with the paper \cite{Ni14} on three-locally-symmetric spaces which were called generalized Wallach spaces. The classification of three-locally-symmetric spaces was obtained in \cite{Ni14} based on the classification of $\mathbb{Z}_2\times \mathbb{Z}_2$-symmetric spaces \cite{BG1,HY1,Ko1}.
\end{remark}

\section{Three-locally-symmetric spaces $G/H$ and involution pairs of $G$}\label{tlss-ip}

Assume that $G/H$ is a three-locally-symmetric space. Then
$$\mathfrak g=\mathfrak h\oplus\mathfrak p_1\oplus \mathfrak p_2\oplus\mathfrak p_3,$$
and it is easy to see that the Lie brackets satisfy
\begin{equation}\label{bracked}
[\mathfrak h,\mathfrak p_i]\subset \mathfrak p_i,\quad [\mathfrak p_i,\mathfrak p_i]\subset \mathfrak h,\quad [\mathfrak p_i,\mathfrak p_j]\subset \mathfrak p_k
\end{equation}
for any $i\in \{1,2,3\}$ and $\{i,j,k\}=\{1,2,3\}$. Define a linear map $\theta_1$ on $\mathfrak g$ by
$$\theta_1|_{\mathfrak h\oplus\mathfrak p_1}=id,\quad \theta_1|_{\mathfrak p_2\oplus\mathfrak p_3}=-id.$$
By the equation~(\ref{bracked}), it is easy to check that $$\theta_1[X,Y]=[\theta_1(X),\theta_1(Y)], \text{ for any } X,Y\in \mathfrak g.$$ It follows that $\theta_1$ is an automorphism, and then an involution of $\mathfrak g$. Similarly, define another linear map $\theta_2$ on $\mathfrak g$ by
$$\theta_2|_{\mathfrak h\oplus\mathfrak p_2}=id,\quad \theta_2|_{\mathfrak p_1\oplus\mathfrak p_3}=-id,$$ which is also an involution of $\mathfrak g$. Moreover we have $\theta_1\theta_2=\theta_2\theta_1$, ${\mathfrak h}=\{X\in {\mathfrak g}|\theta_1(X)=X, \theta_2(X)=X\}$, ${\mathfrak p_1}=\{X\in {\mathfrak g}|\theta_1(X)=X, \theta_2(X)=-X\}$, ${\mathfrak p_2}=\{X\in {\mathfrak g}|\theta_1(X)=-X, \theta_2(X)=X\}$, and ${\mathfrak p_3}=\{X\in {\mathfrak g}|\theta_1(X)=-X, \theta_2(X)=-X\}$.

On the other hand, let $G$ be a compact semisimple connected Lie group with the Lie algebra $\mathfrak g$ and $\rho,\varphi$ be involutions of $\mathfrak g$ satisfying $\rho\varphi=\varphi\rho$. Then we have a decomposition of $\mathfrak g$
$${\mathfrak g}={\mathfrak h}\oplus{\mathfrak p_1}\oplus{\mathfrak p_2}\oplus{\mathfrak p_3},$$
corresponding to $\rho,\varphi$, where ${\mathfrak h}=\{X\in {\mathfrak g}|\rho(X)=X, \varphi(X)=X\}$, ${\mathfrak p_1}=\{X\in {\mathfrak
g}|\rho(X)=X, \varphi(X)=-X\}$, ${\mathfrak p_2}=\{X\in {\mathfrak g}|\rho(X)=-X, \varphi(X)=X\}$, and ${\mathfrak p_3}=\{X\in {\mathfrak g}|\rho(X)=-X, \varphi(X)=-X\}$. It is easy to check that
$$[\mathfrak h,\mathfrak p_i]\subset \mathfrak p_i,\quad [\mathfrak p_i,\mathfrak p_i]\subset \mathfrak h,\quad [\mathfrak p_i,\mathfrak p_j]\subset \mathfrak p_k
$$
for any $i\in \{1,2,3\}$ and $\{i,j,k\}=\{1,2,3\}$. Let $H$ denote the connected Lie subgroup of $G$ with the Lie algebra $\mathfrak h$. If every $\mathfrak p_i$ for $i\in \{1,2,3\}$ is an irreducible $\rm{ad}\mathfrak h$-module, then $G/H$ is a three-locally-symmetric space.

In summary, there is a one-to-one correspondence between the set of three-locally-symmetric spaces and the set of commutative involution pairs of $\mathfrak g$ such that every $\mathfrak p_i$ for $i\in \{1,2,3\}$ is an irreducible $\rm{ad}\mathfrak h$-module.

\section{The classification of three-locally-symmetric spaces}\label{cla}
The section is to give the classification of three-locally-symmetric spaces for a compact simple Lie group. By the discussion in Section~\ref{tlss-ip}, it turns to the classification of certain commutative involution pairs.

Let $G$ be a compact simple connected Lie group with the Lie algebra
${\mathfrak g}$ and $(\theta, \tau)$ be an involution pair of $G$ with $\theta\tau=\tau\theta$. Then for $\theta$, we have a
decomposition,
$${\mathfrak g}={\mathfrak k}+{\mathfrak m},$$ where ${\mathfrak
k}=\{X\in {\mathfrak g}|\theta(X)=X\}$ and ${\mathfrak m}=\{X\in
{\mathfrak g}|\theta(X)=-X\}$. Since $\theta\tau=\tau\theta$, we have that $$\tau(X)\in {\mathfrak k} \text{ for any } X\in {\mathfrak k},$$ which implies that $\tau|_{\mathfrak k}$ is an involution of $\mathfrak k$. Roughly to say, we can give the classification of commutative involution pairs of $\mathfrak g$ by studying the extension of an involution of $\mathfrak k$ to $\mathfrak g$. But an important problem is when an involution of $\mathfrak k$ can be extended to an involution of $\mathfrak g$.

Cartan and Gantmacher made great attributions on the classification of involutions on compact Lie groups. The theory on the extension of involutions of $\mathfrak k$ to $\mathfrak g$ can be found in \cite{Be2}, which is different in the method from that in \cite{Ya1}. There are also some related discussion in \cite{CL,CL1,CH}. The following are the theories without proof.

Let ${\mathfrak t_1}$ be a Cartan
subalgebra of ${\mathfrak k}$ and let $\mathfrak t$ be a Cartan subalgebra of
${\mathfrak g}$ containing $\mathfrak t_1$.

\begin{theorem}[Gantmacher Theorem]\label{chooseH} Let the notations be as above.
Then $\theta$ is conjugate with $\theta_0e^{\mathrm{ad}H}$ under
$\mathrm{Aut} {\mathfrak g}$, where $\theta_0$ is a canonique involution and $H\in {\mathfrak t_1}$.
\end{theorem}

Let $\Pi=\{\alpha_1,\cdots,\alpha_n\}$ be a fundamental system of
${\mathfrak t}$ and $\phi=\sum_{i=1}^nm_i\alpha_i$ be the maximal
root respectively. Let $\alpha_i'=\frac{1}{2}(\alpha_i+\theta_0(\alpha_i))$. Then $\Pi'=\{\alpha_1',\cdots,\alpha_l'\}$ consisting different elements in $\{\alpha_1',\cdots,\alpha_n'\}$ is a fundamental system of $\mathfrak g_0$, where $\mathfrak g_0=\{X\in {\mathfrak g}|\theta_0(X)=X\}$. Denote by $\phi'=\sum_{i=1}^l m_i'\alpha_i'$ the maximal root of $\mathfrak g_0$ respectively. Furthermore we have
\begin{theorem}[\cite{Ya1}]\label{choosealpha} If $H\not=0$, then for some $i$, we can take $H$ satisfying
\begin{equation}\label{H}
\alpha_i'=\alpha_i; \quad \langle H,\alpha_i'\rangle=\pi\sqrt{-1};\quad \langle
H,\alpha_j'\rangle=0, \forall j\not=i.\end{equation} Here $m_i'=1$ or
$m_i'=2$.
\end{theorem}

Moreover, $\mathfrak k$ is described as follows.

\begin{theorem}[\cite{Ya1}]\label{subalgebra}
Let the notations be as above. Assume that $\alpha_i$ satisfies the
identity~(\ref{H}).
\begin{enumerate}
  \item If $\theta_0=Id$ and $m_i=1$, then $\Pi-\{\alpha_i\}$ is the fundamental
  system of $\mathfrak k$, and $\phi$ and $-\alpha_i$ are the
  highest weights of $\mathrm{ad}_{\mathfrak m^{\mathbb C}}{\mathfrak k}$.
  \item If $\theta_0=Id$ and $m_i=2$, then $\Pi-\{\alpha_i\}\cup\{-\phi\}$ is the fundamental
  system of $\mathfrak k$, and $-\alpha_i$ is the
  highest weight of $\mathrm{ad}_{\mathfrak m^{\mathbb C}}{\mathfrak k}$.
  \item If $\theta_0\not=Id$, then $\Pi'-\{\alpha_i'\}\cup\{\beta_0\}$ is the fundamental
  system of $\mathfrak k$, and $-\alpha_i$ is the
  highest weight of $\mathrm{ad}_{\mathfrak m^{\mathbb C}}{\mathfrak k}$.
\end{enumerate}
\end{theorem}

\begin{remark}
In Theorem~\ref{subalgebra}, the dimension of $C(\mathfrak k)$, i.e., the center of $\mathfrak k$, is 1 for case (1); 0 for cases (2) and (3), $\beta_0$ in case (3) is the highest weight of $\mathrm{ad}_{\mathfrak m^{\mathbb C}}{\mathfrak k}$ for $\theta=\theta_0$ corresponding to $\Pi'$.
\end{remark}

Now for any involution $\tau^{\mathfrak k}$ of $\mathfrak k$, we can write
$$\tau^{\mathfrak k}=\tau^{\mathfrak k}_0e^{\mathrm{ad}H^{\mathfrak k}},$$ where $\tau^{\mathfrak k}_0$
is a canonique involution on $\mathfrak k$, $H^{\mathfrak k}\in{\mathfrak t_1}$ and $\tau^{\mathfrak k}_0(H^{\mathfrak k})=H^{\mathfrak k}$.
Since $e^{\mathrm{ad}H^{\mathfrak k}}$ is an inner-automorphism, naturally we can extend $e^{\mathrm{ad}H^{\mathfrak k}}$ to be an automorphism of $\mathfrak g$. Moreover,
\begin{theorem}[\cite{Ya1}]
The canonique involution $\tau^{\mathfrak k}_0$ can be extended to be an automorphism of $\mathfrak g$ if and only if $\tau^{\mathfrak k}_0$ keeps the weight system of $\mathrm{ad}_{\mathfrak m^{\mathbb C}}{\mathfrak k}$ invariant.
\end{theorem}

If $C(\mathfrak k)\not=0$, then $\dim C(\mathfrak k)=1$. Thus $\tau^{\mathfrak k}_0(Z)=Z$ or $\tau^{\mathfrak k}_0(Z)=-Z$ for any $Z\in C(\mathfrak k)$.

\begin{theorem}[\cite{Ya1}]
Assume that $C(\mathfrak k)\not=0$ and $\tau^{\mathfrak k}_0(Z)=Z$ for any $Z\in C(\mathfrak k)$. If $\tau^{\mathfrak k}$ can be extended to be an automorphism of $\mathfrak g$, then $\tau^{\mathfrak k}$ can be extended to be an involution of $\mathfrak g$.
\end{theorem}

For the other cases, we have the following theorems.
\begin{theorem}[\cite{Ya1}]
Assume that $C(\mathfrak k)=0$, or $C(\mathfrak k)\not=0$ but $\tau^{\mathfrak k}_0(Z)=-Z$ for any $Z\in C(\mathfrak k)$. If $\tau$ is an automorphism of $\mathfrak g$ extended by an involution $\tau^{\mathfrak k}$ of $\mathfrak k$, then $\tau^2=Id$ or $\tau^2=\theta$. Furthermore, the following conditions are equivalent:
\begin{enumerate}
   \item There exists an automorphism $\tau$ of $\mathfrak g$ extended by $\tau^{\mathfrak k}$ which is an involution.
   \item Every automorphism $\tau$ of $\mathfrak g$ extended by $\tau^{\mathfrak k}$ is an involution.
\end{enumerate}
\end{theorem}

Then it is enough to determine when the automorphism extended by $\tau^{\mathfrak k}$ is an involution.
\begin{theorem}[\cite{Ya1}]
Let $\tau_0$ be the automorphism of $\mathfrak g$ extended by the canonique involution $\tau^{\mathfrak k}_0$ on $\mathfrak k$. Then $\tau_0^2=Id$ except $\mathfrak g=A_n^i$ and $n$ is even. For $e^{\mathrm{ad}H^{\mathfrak k}}$, we have:
\begin{enumerate}
 \item If $\theta_0\not=Id$, then the natural extension of $e^{\mathrm{ad}H^{\mathfrak k}}$ is an involution.
 \item Assume that $\theta_0=Id$. Let $\alpha'_{i_1},\cdots,\alpha'_{i_k}$ be the roots satisfying $\langle \alpha'_j,H\rangle\not=0$. Then the natural extension of $e^{\mathrm{ad}H^{\mathfrak k}}$ is an involution if and only if $\sum_{j=i_1}^{i_k}m'_j$ is even.
\end{enumerate}
\end{theorem}

In particular, for every cases in Theorem~\ref{subalgebra},
\begin{theorem}[\cite{Be2,Ya1}]\label{lem}
If $\tau$ is an involution of $\mathfrak g$ extended by an involution $\tau^{\mathfrak k}$ on $\mathfrak k$, then every extension of $\tau^{\mathfrak k}$ is an involution of $\mathfrak k$, which is equivalent with $\tau$ or $\tau\theta$.
\end{theorem}

Up to now, we can obtain the classification of commutative involution pairs by that of $(\theta,\tau^{\mathfrak k})$ based on the above theory. By Theorem~\ref{lem}, for an involution $\tau^{\mathfrak k}$ on $\mathfrak k$ which can be extended to an involution $\tau$ of $\mathfrak g$, we have two involution pairs $(\theta,\tau)$ and $(\theta,\theta\tau)$ which determine the same three-locally-symmetric space. So, without loss of generality, we denote by $\tau$ the natural extension of $\tau^{\mathfrak k}$. Let ${\mathfrak h}=\{X\in {\mathfrak k}|\tau^{\mathfrak k}(X)=X\}$, ${\mathfrak p}_1=\{X\in {\mathfrak k}|\tau(X)=-X\}$, ${\mathfrak p}_2=\{X\in {\mathfrak m}|\tau(X)=X\}$, and ${\mathfrak p}_3=\{X\in {\mathfrak m}|\tau(X)=-X\}$. We shall pick up certain pairs $(\theta,\tau^{\mathfrak k})$ by the following steps.

{\bf Step 1}:  we obtain the classification of the pair $(\theta,\tau^{\mathfrak k})$ satisfying that the extension $\tau$ of $\tau^{\mathfrak k}$ is an involution and ${\mathfrak p}_1,{\mathfrak p}_2,{\mathfrak p}_3$ are irreducible as $\mathrm{ad}{\mathfrak h}$-modules.

{\bf Step 2}: among the pairs given by Step 1, we remain only one if several pairs determine the same three-locally-symmetric space.

By the above steps, we obtain the list of the following pairs $(\theta,\tau^{\mathfrak k})$. Furthermore, we make some remarks on them.

\medskip
Case $G=A_l(l\geq1)$. The Dynkin diagram with the maximal root is

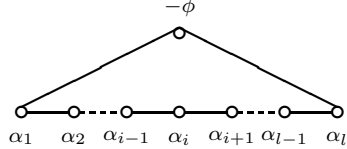
\begin{center}
\setlength{\unitlength}{0.7mm}
\begin{picture}(70,30)(0,0)
\thicklines \multiput(0,10)(10,0){7}{\circle{2}}
\multiput(1,10)(10,0){1}{\line(1,0){8}}
\multiput(11,10)(2,0){4}{\line(1,0){1}}
\multiput(21,10)(10,0){2}{\line(1,0){8}}
\multiput(41,10)(2,0){4}{\line(1,0){1}}
\multiput(51,10)(10,0){1}{\line(1,0){8}} \put(30,25){\circle{2}}
\put(0,11){\line(2,1){30}} \put(60,11){\line(-2,1){30}}
\put(0,5){\makebox(0,0){$\scriptstyle {\alpha_1}$}}
\put(10,5){\makebox(0,0){$\scriptstyle {\alpha_2}$}}
\put(20,5){\makebox(0,0){$\scriptstyle {\alpha_{i-1}}$}}
\put(30,5){\makebox(0,0){$\scriptstyle {\alpha_{i}}$}}
\put(40,5){\makebox(0,0){$\scriptstyle {\alpha_{i+1}}$}}
\put(50,5){\makebox(0,0){$\scriptstyle {\alpha_{l-1}}$}}
\put(60,5){\makebox(0,0){$\scriptstyle {\alpha_{l}}$}}
\put(30,30){\makebox(0,0){$\scriptstyle {-\phi}$}}
\end{picture}
\end{center}

$(InP1)$ $l=1$, $\theta=e^{\mathrm{ad}H}$, where $\langle
H,\alpha_1\rangle=\pi\sqrt{-1}$; $\tau^{\mathfrak k}|_{C(\mathfrak k)}=-Id$.

$(InP2)$ $l\geq 3$ is odd, $\theta=e^{\mathrm{ad}H}$, where $\langle
H,\alpha_{\frac{l+1}{2}}\rangle=\pi\sqrt{-1}$ and
$\langle H,\alpha_k\rangle=0$ for any $\alpha_k\in
\Pi-\{\alpha_i\}$; $\tau^{\mathfrak k}(\alpha_k)=\alpha_{l+1-k}$ for $k\not=\frac{l+1}{2}$ and $\tau^{\mathfrak k}|_{C(\mathfrak k)}=Id$.

$(InP3)$ $l\geq 2$, $\theta=e^{\mathrm{ad}H}$, where $\langle
H,\alpha_i\rangle=\pi\sqrt{-1}$ for some $1\leq i\leq [\frac{l+1}{2}]$ and
$\langle H,\alpha_k\rangle=0$ for any $\alpha_k\in
\Pi-\{\alpha_i\}$; $\tau^{\mathfrak k}=e^{\mathrm{ad}H_1}$, where
$\langle H_1,\alpha_j\rangle=\pi\sqrt{-1}$ for some $\alpha_j\in
\Pi-\{\alpha_i\}$ and $\langle H,\alpha_k\rangle=0$ for any
$\alpha_k\in \Pi-\{\alpha_i,\alpha_j\}$. Furthermore, we may require that $1\leq i\leq [\frac{l+1}{3}]$ and $2i\leq j\leq [\frac{l+i+1}{2}]$ for Step 2.

\medskip
Case $G=B_l(l\geq2)$. The Dynkin diagram with the maximal root is
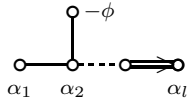
\begin{center}
\setlength{\unitlength}{0.7mm}
\begin{picture}(50,23)(0,0) \thicklines
\multiput(0,10)(10,0){4}{\circle{2}} \put(10,20){\circle{2}}
\put(30,10){\circle{2}} \multiput(1,10)(20,0){1}{\line(1,0){8}}
\multiput(11,10)(2,0){4}{\line(1,0){1}}
\multiput(21,9.5)(0,1){2}{\line(1,0){8}} \put(10,11){\line(0,1){8}}
\put(25.5,8.5){$>$} \put(0,5){\makebox(0,0){$\scriptstyle
\alpha_1$}} \put(10,5){\makebox(0,0){$\scriptstyle \alpha_2$}}
\put(15,20){\makebox(0,0){$\scriptstyle -\phi$}}
\put(30,5){\makebox(0,0){$\scriptstyle \alpha_l$}}
\end{picture}
\end{center}

$(InP4)$ $\theta=e^{\mathrm{ad}H}$, where $\langle
H,\alpha_i\rangle=\pi\sqrt{-1}$ for some $2\leq i\leq l$ and
$\langle H,\alpha_k\rangle=0$ for any $\alpha_k\in
\Pi-\{\alpha_i\}$; $\tau^{\mathfrak k}=e^{\mathrm{ad}H_1}$, where
$\langle H_1,\alpha_j\rangle=\pi\sqrt{-1}$ for some $2\leq j< i$ and $\langle H,\alpha_k\rangle=0$ for any
$\alpha_k\in \Pi-\{\alpha_i,\alpha_j\}\cup{\{-\phi\}}$. Moreover, we may assume that $2< i\leq l$ and $\frac{i}{2}\leq j\leq i-1$.

$(InP5)$ $\theta=e^{\mathrm{ad}H}$, where $\langle
H,\alpha_i\rangle=\pi\sqrt{-1}$ for some $2\leq i\leq l$ and
$\langle H,\alpha_k\rangle=0$ for any $\alpha_k\in
\Pi-\{\alpha_i\}$; $\tau^{\mathfrak k}(\alpha_1)=-\phi$, $\tau^{\mathfrak k}(-\phi)=\alpha_1$, and $\tau^{\mathfrak k}(\alpha_k)=\alpha_k$ for any
$\alpha_k\in \Pi-\{\alpha_1,\alpha_i\}$. Moreover, we may assume that $\frac{l+1}{2}\leq i\leq l$.

$(InP6)$ $\theta=e^{\mathrm{ad}H}$, where $\langle
H,\alpha_i\rangle=\pi\sqrt{-1}$ for some $2\leq i\leq l$ and
$\langle H,\alpha_k\rangle=0$ for any $\alpha_k\in
\Pi-\{\alpha_i\}$; $\tau^{\mathfrak k}=\tau^{\mathfrak k}_0e^{\mathrm{ad}H_1}$, where $\tau^{\mathfrak k}_0(\alpha_1)=-\phi$, $\tau^{\mathfrak k}_0(-\phi)=\alpha_1$, and $\tau^{\mathfrak k}_0(\alpha_k)=\alpha_k$ for any
$\alpha_k\in \Pi-\{\alpha_1,\alpha_i\}$, $\langle H_1,\alpha_j\rangle=\pi\sqrt{-1}$ for some $2\leq j<i$ and $\langle H,\alpha_k\rangle=0$ for any
$\alpha_k\in \Pi-\{\alpha_i,\alpha_j\}\cup{\{-\phi\}}$. Furthermore, we may require that $[\frac{2l+3}{3}]\leq i\leq l$ and $[\frac{i+2}{2}]\leq j\leq 2i-l$.

\medskip
Case $G=C_l(l\geq3)$. The Dynkin diagram with the maximal root is
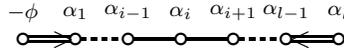
\begin{center}
\setlength{\unitlength}{0.7mm}
\begin{picture}(60,15)(0,-2)
\thicklines \multiput(0,0)(10,0){7}{\circle{2}}
\multiput(1,-0.5)(0,1){2}{\line(1,0){8}}
\multiput(11,0)(2,0){4}{\line(1,0){1}}
\multiput(21,0)(10,0){2}{\line(1,0){8}}
\multiput(41,0)(2,0){4}{\line(1,0){1}}
\multiput(51,-0.5)(0,1){2}{\line(1,0){8}} \put(50.5,-1.5){$<$}
\put(5.5,-1.5){$>$} \put(0,5){\makebox(0,0){$\scriptstyle -\phi$}}
\put(10,5){\makebox(0,0){$\scriptstyle \alpha_1$}}
\put(20,5){\makebox(0,0){$\scriptstyle \alpha_{i-1}$}}
\put(30,5){\makebox(0,0){$\scriptstyle \alpha_i$}}
\put(40,5){\makebox(0,0){$\scriptstyle \alpha_{i+1}$}}
\put(50,5){\makebox(0,0){$\scriptstyle \alpha_{l-1}$}}
\put(60,5){\makebox(0,0){$\scriptstyle \alpha_l$}}
\end{picture}
\end{center}

$(InP7)$ $\theta=e^{\mathrm{ad}H}$, where $\langle
H,\alpha_i\rangle=\pi\sqrt{-1}$ for some $1\leq i\leq [\frac{l}{2}]$ and
$\langle H,\alpha_k\rangle=0$ for any $\alpha_k\in
\Pi-\{\alpha_i\}$; $\tau^{\mathfrak k}=e^{\mathrm{ad}H_1}$, where
$\langle H_1,\alpha_j\rangle=\pi\sqrt{-1}$ for some $\alpha_j\in
\Pi-\{\alpha_i,\alpha_l,-\phi\}$ and $\langle H,\alpha_k\rangle=0$ for any
$\alpha_k\in \Pi-\{\alpha_i,\alpha_j\}\cup{\{-\phi\}}$. Furthermore, we may require that $1\leq i\leq [\frac{l}{3}]$ and $2i\leq j\leq [\frac{l+i}{2}]$.

\medskip
Case $G=D_l(l\geq4)$. The Dynkin diagram with the maximal root is
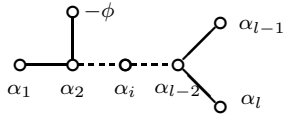
\begin{center}
\setlength{\unitlength}{0.7mm}
\begin{picture}(50,25)(0,-8) \thicklines
\multiput(0,0)(10,0){4}{\circle{2}}
\multiput(38,-8)(0,16){2}{\circle{2}} \put(10,10){\circle{2}}
\put(20,0){\circle{2}} \multiput(1,0)(20,0){1}{\line(1,0){8}}
\multiput(11,0)(2,0){4}{\line(1,0){1}}
\multiput(21,0)(2,0){4}{\line(1,0){1}} \put(31,1){\line(1,1){6}}
\put(31,-1){\line(1,-1){6}} \put(10,1){\line(0,1){8}}
\put(0,-5){\makebox(0,0){$\scriptstyle \alpha_1$}}
\put(10,-5){\makebox(0,0){$\scriptstyle \alpha_2$}}
\put(15,10){\makebox(0,0){$\scriptstyle -\phi$}}
\put(20,-5){\makebox(0,0){$\scriptstyle \alpha_i$}}
\put(30,-5){\makebox(0,0){$\scriptstyle \alpha_{l-2}$}}
\put(46,7){\makebox(0,0){$\scriptstyle \alpha_{l-1}$}}
\put(44,-7){\makebox(0,0){$\scriptstyle \alpha_l$}}
\end{picture}
\end{center}

$(InP8)$ $\theta=e^{\mathrm{ad}H}$, where $\langle
H,\alpha_i\rangle=\pi\sqrt{-1}$ for some $2\leq i\leq [\frac{l}{2}]$ and
$\langle H,\alpha_k\rangle=0$ for any $\alpha_k\in
\Pi-\{\alpha_i\}$; $\tau^{\mathfrak k}=e^{\mathrm{ad}H_1}$, where
$\langle H_1,\alpha_j\rangle=\pi\sqrt{-1}$ for some $\alpha_j\in
\Pi-\{\alpha_1,\alpha_i,\alpha_{l-1},\alpha_l\}$ and $\langle
H,\alpha_k\rangle=0$ for any $\alpha_k\in
\Pi-\{\alpha_i,\alpha_j\}\cup{\{-\phi\}}$. Furthermore, we may require that $1\leq i\leq [\frac{l}{3}]$ and $2i\leq j\leq [\frac{l+i}{2}]$

$(InP9)$ $\theta=e^{\mathrm{ad}H}$, where $\langle
H,\alpha_i\rangle=\pi\sqrt{-1}$ for some $1\leq i\leq l-2$ and
$\langle H,\alpha_k\rangle=0$ for any $\alpha_k\in
\Pi-\{\alpha_i\}$; $\tau^{\mathfrak k}(\alpha_l)=\alpha_{l-1}$ and $\tau^{\mathfrak k}(\alpha_{l-1})=\alpha_{l}$.

$(InP10)$ $\theta=e^{\mathrm{ad}H}$, where $\langle
H,\alpha_i\rangle=\pi\sqrt{-1}$ for some $1\leq i\leq l-2$ and
$\langle H,\alpha_k\rangle=0$ for any $\alpha_k\in
\Pi-\{\alpha_i\}$; $\tau^{\mathfrak k}=\tau^{\mathfrak k}_0e^{\mathrm{ad}H_1}$, where $\tau^{\mathfrak k}_0(\alpha_l)=\alpha_{l-1}$ and $\tau^{\mathfrak k}_0(\alpha_{l-1})=\alpha_{l}$, $\langle H_1,\alpha_j\rangle=\pi\sqrt{-1}$ for some $i<j\leq \frac{l+i-1}{2}$ and $\langle
H,\alpha_k\rangle=0$ for any $\alpha_k\in
\Pi-\{\alpha_i,\alpha_j\}\cup{\{-\phi\}}$.

$(InP11)$ $\theta=e^{\mathrm{ad}H}$, where $\langle
H,\alpha_1\rangle=\pi\sqrt{-1}$ and $\langle H,\alpha_k\rangle=0$
for any $\alpha_k\in \Pi-\{\alpha_1\}$; $\tau^{\mathfrak k}|_{C(\mathfrak k)}=-Id$.

$(InP12)$ $\theta=e^{\mathrm{ad}H}$, where $\langle
H,\alpha_1\rangle=\pi\sqrt{-1}$ and $\langle H,\alpha_k\rangle=0$
for any $\alpha_k\in \Pi-\{\alpha_1\}$; $\tau^{\mathfrak
k}=e^{\mathrm{ad}H_1}$, where $\langle
H_1,\alpha_l\rangle=\pi\sqrt{-1}$ and $\langle H,\alpha_k\rangle=0$
for any $\alpha_k\in \Pi-\{\alpha_1,\alpha_l\}$.

\medskip
Case $G=E_6$. The Dynkin diagram with the maximal root is
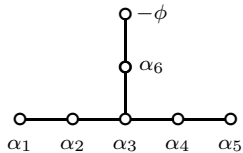
\begin{center}
\setlength{\unitlength}{0.7mm}
\begin{picture}(50,33)(0,0)
\thicklines \multiput(0,10)(10,0){5}{\circle{2}}
\multiput(20,20)(0,10){2}{\circle{2}} \put(20,20){\circle{2}}
\multiput(1,10)(10,0){4}{\line(1,0){8}}
\multiput(20,11)(0,10){2}{\line(0,1){8}}
\put(0,5){\makebox(0,0){$\scriptstyle \alpha_1$}}
\put(25,20){\makebox(0,0){$\scriptstyle \alpha_6$}}
\put(10,5){\makebox(0,0){$\scriptstyle \alpha_2$}}
\put(20,5){\makebox(0,0){$\scriptstyle \alpha_3$}}
\put(30,5){\makebox(0,0){$\scriptstyle \alpha_4$}}
\put(40,5){\makebox(0,0){$\scriptstyle \alpha_5$}}
\put(25,30){\makebox(0,0){$\scriptstyle -\phi$}}
\end{picture}
\end{center}

$(InP13)$ $\theta=e^{\mathrm{ad}H}$, where $\langle
H,\alpha_1\rangle=\pi\sqrt{-1}$ and $\langle H,\alpha_k\rangle=0$
for any $\alpha_k\in \Pi-\{\alpha_1\}$; $\tau^{\mathfrak
k}=e^{\mathrm{ad}H_1}$, where $\langle
H_1,\alpha_5\rangle=\pi\sqrt{-1}$ and $\langle H,\alpha_k\rangle=0$
for any $\alpha_k\in \Pi-\{\alpha_1,\alpha_5\}$.

$(InP14)$ $\theta=e^{\mathrm{ad}H}$, where $\langle
H,\alpha_6\rangle=\pi\sqrt{-1}$ and $\langle H,\alpha_k\rangle=0$
for any $\alpha_k\in \Pi-\{\alpha_6\}$; $\tau^{\mathfrak
k}=e^{\mathrm{ad}H_1}$, where $\langle
H_1,\alpha_2\rangle=\pi\sqrt{-1}$ and $\langle H,\alpha_k\rangle=0$
for any $\alpha_k\in \Pi-\{\alpha_2,\alpha_6\}\cup\{-\phi\}$.

$(InP15)$ $\theta=e^{\mathrm{ad}H}$, where $\langle
H,\alpha_6\rangle=\pi\sqrt{-1}$ and $\langle H,\alpha_k\rangle=0$
for any $\alpha_k\in \Pi-\{\alpha_6\}$; $\tau^{\mathfrak
k}(\alpha_i)=\alpha_{6-i}$ for $i=1,2,3,4,5$.

\medskip
Case $G=E_7$. The Dynkin diagram with the maximal root is
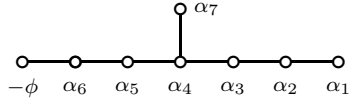
\begin{center}
\setlength{\unitlength}{0.7mm}
\begin{picture}(65,25)(0,0)
\thicklines \multiput(0,10)(10,0){7}{\circle{2}}
\put(30,20){\circle{2}} \put(10,10){\circle{2}}
\multiput(1,10)(10,0){6}{\line(1,0){8}} \put(30,11){\line(0,1){8}}
\put(10,5){\makebox(0,0){$\scriptstyle \alpha_6$}}
\put(35,20){\makebox(0,0){$\scriptstyle \alpha_7$}}
\put(20,5){\makebox(0,0){$\scriptstyle \alpha_5$}}
\put(30,5){\makebox(0,0){$\scriptstyle \alpha_4$}}
\put(40,5){\makebox(0,0){$\scriptstyle \alpha_3$}}
\put(50,5){\makebox(0,0){$\scriptstyle \alpha_2$}}
\put(60,5){\makebox(0,0){$\scriptstyle \alpha_1$}}
\put(0,5){\makebox(0,0){$\scriptstyle -\phi$}}
\end{picture}
\end{center}

$(InP16)$ $\theta=e^{\mathrm{ad}H}$, where $\langle
H,\alpha_6\rangle=\pi\sqrt{-1}$ and $\langle H,\alpha_k\rangle=0$
for any $\alpha_k\in \Pi-\{\alpha_6\}$; $\tau^{\mathfrak
k}=e^{\mathrm{ad}H_1}$, where $\langle
H_1,\alpha_2\rangle=\pi\sqrt{-1}$ and $\langle H,\alpha_k\rangle=0$
for any $\alpha_k\in \Pi-\{\alpha_2,\alpha_6\}\cup\{-\phi\}$.

$(InP17)$ $\theta=e^{\mathrm{ad}H}$, where $\langle
H,\alpha_7\rangle=\pi\sqrt{-1}$ and $\langle H,\alpha_k\rangle=0$
for any $\alpha_k\in \Pi-\{\alpha_7\}$; $\tau^{\mathfrak
k}=e^{\mathrm{ad}H_1}$, where $\langle
H_1,\alpha_2\rangle=\pi\sqrt{-1}$ and $\langle H,\alpha_k\rangle=0$
for any $\alpha_k\in \Pi-\{\alpha_2,\alpha_7\}\cup\{-\phi\}$.

$(InP18)$ $\theta=e^{\mathrm{ad}H}$, where $\langle
H,\alpha_7\rangle=\pi\sqrt{-1}$ and $\langle H,\alpha_k\rangle=0$
for any $\alpha_k\in \Pi-\{\alpha_7\}$; $\tau^{\mathfrak k}=\tau^{\mathfrak k}_0e^{\mathrm{ad}H_1}$, where $\tau^{\mathfrak
k}_0(\alpha_i)=\alpha_{8-i}$ for $i=2,3,4,5,6$, $\tau^{\mathfrak
k}_0(\alpha_1)=-\phi$, and $\tau^{\mathfrak
k}_0(-\phi)=\alpha_1$, $\langle H_1,\alpha_4\rangle=\pi\sqrt{-1}$ and $i<j\leq l-3$ and $\langle
H_1,\alpha_k\rangle=0$ for any $\alpha_k\in
\Pi-\{\alpha_4,\alpha_7\}\cup\{-\phi\}$.

\medskip
Case $G=E_8$. The Dynkin diagram with the maximal root is
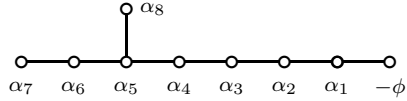
\begin{center}
\setlength{\unitlength}{0.7mm}
\begin{picture}(75,25)(0,0)
\thicklines \multiput(0,10)(10,0){8}{\circle{2}}
\put(20,20){\circle{2}} \multiput(1,10)(10,0){7}{\line(1,0){8}}
\put(20,11){\line(0,1){8}} \put(60,10){\circle{2}}
\put(0,5){\makebox(0,0){$\scriptstyle \alpha_7$}}
\put(25,20){\makebox(0,0){$\scriptstyle \alpha_8$}}
\put(10,5){\makebox(0,0){$\scriptstyle \alpha_6$}}
\put(20,5){\makebox(0,0){$\scriptstyle \alpha_5$}}
\put(30,5){\makebox(0,0){$\scriptstyle \alpha_4$}}
\put(40,5){\makebox(0,0){$\scriptstyle \alpha_3$}}
\put(50,5){\makebox(0,0){$\scriptstyle \alpha_2$}}
\put(60,5){\makebox(0,0){$\scriptstyle \alpha_1$}}
\put(70,5){\makebox(0,0){$\scriptstyle -\phi$}}
\end{picture}
\end{center}

$(InP19)$ $\theta=e^{\mathrm{ad}H}$, where $\langle
H,\alpha_7\rangle=\pi\sqrt{-1}$ and $\langle H,\alpha_k\rangle=0$
for any $\alpha_k\in \Pi-\{\alpha_7\}$; $\tau^{\mathfrak
k}=e^{\mathrm{ad}H_1}$, where $\langle
H_1,\alpha_1\rangle=\pi\sqrt{-1}$ and $\langle H,\alpha_k\rangle=0$
for any $\alpha_k\in \Pi-\{\alpha_1,\alpha_7\}\cup\{-\phi\}$.

$(InP20)$ $\theta=e^{\mathrm{ad}H}$, where $\langle
H,\alpha_7\rangle=\pi\sqrt{-1}$ and $\langle H,\alpha_k\rangle=0$
for any $\alpha_k\in \Pi-\{\alpha_7\}$; $\tau^{\mathfrak
k}=e^{\mathrm{ad}H_1}$, where $\langle
H_1,\alpha_3\rangle=\pi\sqrt{-1}$ and $\langle H,\alpha_k\rangle=0$
for any $\alpha_k\in \Pi-\{\alpha_3,\alpha_7\}\cup\{-\phi\}$.
\medskip

Case $G=F_4$. The Dynkin diagram with the maximal root is
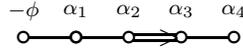
\begin{center}
\setlength{\unitlength}{0.7mm}
\begin{picture}(40,10)(0,0)
\thicklines \multiput(0,0)(10,0){5}{\circle{2}}
\put(10,0){\circle{2}} \multiput(1,0)(10,0){2}{\line(1,0){8}}
\multiput(21,-0.5)(0,1){2}{\line(1,0){8}} \put(31,0){\line(1,0){8}}
\put(25.5,-1.5){$>$} \put(0,5){\makebox(0,0){$\scriptstyle -\phi$}}
\put(10,5){\makebox(0,0){$\scriptstyle \alpha_1$}}
\put(20,5){\makebox(0,0){$\scriptstyle \alpha_2$}}
\put(30,5){\makebox(0,0){$\scriptstyle \alpha_3$}}
\put(40,5){\makebox(0,0){$\scriptstyle \alpha_4$}}
\end{picture}
\end{center}

$(InP21)$ $\theta=e^{\mathrm{ad}H}$, where $\langle
H,\alpha_4\rangle=\pi\sqrt{-1}$ and $\langle H,\alpha_k\rangle=0$
for any $\alpha_k\in \Pi-\{\alpha_4\}$; $\tau^{\mathfrak
k}=e^{\mathrm{ad}H_1}$, where $\langle
H_1,\alpha_3\rangle=\pi\sqrt{-1}$ and $\langle H,\alpha_k\rangle=0$
for any $\alpha_k\in \Pi-\{\alpha_3,\alpha_4\}\cup\{-\phi\}$.

$(InP22)$ $\theta=e^{\mathrm{ad}H}$, where $\langle
H,\alpha_4\rangle=\pi\sqrt{-1}$ and $\langle H,\alpha_k\rangle=0$
for any $\alpha_k\in \Pi-\{\alpha_4\}$; $\tau^{\mathfrak
k}=e^{\mathrm{ad}H_1}$, where $\langle
H_1,\alpha_1\rangle=\pi\sqrt{-1}$ and $\langle H,\alpha_k\rangle=0$
for any $\alpha_k\in \Pi-\{\alpha_1,\alpha_4\}\cup\{-\phi\}$.

\begin{remark}
Consider $G=A_l$ when $l\geq 1$ is odd. We take the involution $\theta=e^{\mathrm{ad}H}$, where $\langle
H,\alpha_{\frac{l+1}{2}}\rangle=\pi\sqrt{-1}$ and
$\langle H,\alpha_k\rangle=0$ for any $\alpha_k\in
\Pi-\{\alpha_i\}$. Then $\tau^{\mathfrak k}$ defined by $\tau^{\mathfrak k}(\alpha_k)=\alpha_{\frac{l+1}{2}+k}$ and $\tau^{\mathfrak k}(\alpha_{\frac{l+1}{2}+k})=\alpha_k$ for any $1\leq k\leq \frac{l-1}{2}$ is an involution on $\mathfrak k$. By the above theory, if $\tau^{\mathfrak k}$ can be extended to an involution of $\mathfrak g$, we obtain $\tau^{\mathfrak k}(-\alpha_{\frac{l+1}{2}})=\phi$. It is equivalent to $\tau^{\mathfrak k}|_{C(\mathfrak k)}=-Id$. If $l\geq 3$, $\mathfrak p_1$ is reducible. If $l=1$, we get the case $(InP1)$.
\end{remark}

\begin{remark}
The case $(InP4)$ is valid for $i=1$ and $j$, which determines the same three-locally-symmetric space with $i'=j$ and $j'=j-1$.
\end{remark}

\begin{remark}
If we take $\theta$ the case $(InP4)$ for $i=1$, then $\tau^{\mathfrak k}|_{C(\mathfrak k)}=-Id$ can be extended to an involution of $\mathfrak g$. This pair determines the same three-locally-symmetric space as the case $(InP4)$ for $i=l$.
\end{remark}

\begin{remark}
The involution $\tau^{\mathfrak k}$ in cases $(InP5)$ and $(InP6)$ are outer-automorphism on $\mathfrak k$, which can be extended to an inner-automorphism of $\mathfrak g$.
\end{remark}

\begin{remark}
Consider the following involution pair $(\theta,\tau)$ of $G=C_l$. First for the Dyinkin diagram of $C_l$, $\theta=e^{\mathrm{ad}H}$ is an involution of $\mathfrak g$, where $\langle H,\alpha_l\rangle=\pi\sqrt{-1}$ and others zero. Then $\mathfrak k$ is the direct sum of $A_{l-1}$ and the center of one-dimension. Given an involution $\tau^{\mathfrak k}=e^{\mathrm{ad}H_1}$ of $\mathfrak k$, where
$\langle H_1,\alpha_j\rangle=\pi\sqrt{-1}$ for some $j\in \{1,2,\cdots,l-1\}$ and $\langle H,\alpha_k\rangle=0$ for any
$\alpha_k\in \Pi-\{\alpha_j,\alpha_l\}$. Then the extension of $\tau_{\mathfrak k}$ is an involution of $\mathfrak g$ and
\begin{enumerate}
\item ${\mathfrak p}_1$ is an irreducible $\mathrm{ad}{\mathfrak h}$-module.
\item By a result in \cite{Be2}, the extension of $\tau^{\mathfrak k}$ is $\tau$ or $\tau\theta$. Here $\tau$ is also denoted by $e^{\mathrm{ad}H_1}$, where $\langle H_1,\alpha_j\rangle=\pi\sqrt{-1}$ and others zero.
\end{enumerate}
We can prove that ${\mathfrak p}_2$ is a reducible $\mathrm{ad}{\mathfrak h}$-module. In fact, the Dynkin diagram of $\mathfrak h+\mathfrak p_2$, i.e. the set of fixed points of $\tau$, is
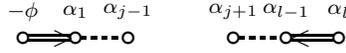
\begin{center}
\setlength{\unitlength}{0.7mm}
\begin{picture}(60,10)(0,-2)
\thicklines \multiput(0,0)(10,0){3}{\circle{2}}
\multiput(40,0)(10,0){3}{\circle{2}}
\multiput(1,-0.5)(0,1){2}{\line(1,0){8}}
\multiput(11,0)(2,0){4}{\line(1,0){1}}
\multiput(41,0)(2,0){4}{\line(1,0){1}}
\multiput(51,-0.5)(0,1){2}{\line(1,0){8}} \put(50.5,-1.5){$<$}
\put(5.5,-1.5){$>$} \put(0,5){\makebox(0,0){$\scriptstyle -\phi$}}
\put(10,5){\makebox(0,0){$\scriptstyle \alpha_1$}}
\put(20,5){\makebox(0,0){$\scriptstyle \alpha_{j-1}$}}
\put(40,5){\makebox(0,0){$\scriptstyle \alpha_{j+1}$}}
\put(50,5){\makebox(0,0){$\scriptstyle \alpha_{l-1}$}}
\put(60,5){\makebox(0,0){$\scriptstyle \alpha_l$}}
\end{picture}
\end{center}
Then by the definition of $\theta$, the Dynkin diagram of $\mathfrak h$ is
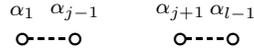
\begin{center}
\setlength{\unitlength}{0.7mm}
\begin{picture}(60,10)(0,-2)
\thicklines \multiput(10,0)(10,0){2}{\circle{2}}
\multiput(40,0)(10,0){2}{\circle{2}}
\multiput(11,0)(2,0){4}{\line(1,0){1}}
\multiput(41,0)(2,0){4}{\line(1,0){1}}
\put(10,5){\makebox(0,0){$\scriptstyle \alpha_1$}}
\put(20,5){\makebox(0,0){$\scriptstyle \alpha_{j-1}$}}
\put(40,5){\makebox(0,0){$\scriptstyle \alpha_{j+1}$}}
\put(50,5){\makebox(0,0){$\scriptstyle \alpha_{l-1}$}}
\end{picture}
\end{center}
and on $\mathfrak h+\mathfrak p_2$, $\langle
H,\alpha_l\rangle=\pi\sqrt{-1}$, $\langle
H,-\phi\rangle=\pi\sqrt{-1}$, and others zero. That is, ${\mathfrak p}_2$ is a reducible $\mathrm{ad}{\mathfrak h}$-module.
\end{remark}

\begin{remark}\label{exam}
This remark is in detail to describe $\mathfrak p_1,\mathfrak p_2,\mathfrak p_3$ as irreducible modules of $\mathrm{ad} \mathfrak h$ for the case $(InP21)$, which is given in \cite{CL1}. For this case, $\mathfrak p_1,\mathfrak p_2,\mathfrak p_3$ are pairwise nonisomorphic. In fact, we can choose a fundamental system
$\{\alpha_1,\alpha_2,\alpha_3,\alpha_4\}$ of $F_4$ such that the
Dynkin diagram of $F_4$ corresponding to the fundamental system is
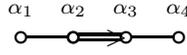
\begin{center}
\setlength{\unitlength}{0.7mm}
\begin{picture}(30,15)(0,-5)
\thicklines \multiput(0,0)(10,0){4}{\circle{2}}
\put(10,0){\circle{2}} \multiput(1,0)(10,0){1}{\line(1,0){8}}
\multiput(11,-0.5)(0,1){2}{\line(1,0){8}} \put(21,0){\line(1,0){8}}
\put(16,-1.5){$>$} \put(0,5){\makebox(0,0){$\scriptstyle \alpha_1$}}
\put(10,5){\makebox(0,0){$\scriptstyle \alpha_2$}}
\put(20,5){\makebox(0,0){$\scriptstyle \alpha_3$}}
\put(30,5){\makebox(0,0){$\scriptstyle \alpha_4$}}
\end{picture}
\end{center}
Consider the involution
$\theta=e^{\mathrm{ad}H}$ defined by $$\langle H,\alpha_4\rangle=\pi\sqrt{-1};\quad \langle H,\alpha_j\rangle=0,
\forall j\not=4.$$ Let
$\phi=2\alpha_1+3\alpha_2+4\alpha_3+2\alpha_4$. Then the Dynkin
diagram of ${\mathfrak k}$ is
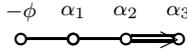
\begin{center}
\setlength{\unitlength}{0.7mm}
\begin{picture}(30,15)(0,-5)
\thicklines \multiput(0,0)(10,0){4}{\circle{2}}
\put(10,0){\circle{2}} \multiput(1,0)(10,0){2}{\line(1,0){8}}
\multiput(21,-0.5)(0,1){2}{\line(1,0){8}} \put(26,-1.5){$>$}
\put(0,5){\makebox(0,0){$\scriptstyle -\phi$}}
\put(10,5){\makebox(0,0){$\scriptstyle \alpha_1$}}
\put(20,5){\makebox(0,0){$\scriptstyle \alpha_2$}}
\put(30,5){\makebox(0,0){$\scriptstyle \alpha_3$}}
\end{picture}
\end{center}
The involution $\tau$ is the extension of $\tau^{\mathfrak k}$,
where $\tau^{\mathfrak k}=e^{\mathrm{ad}H_1}$ satisfies
$$\langle H_1,\alpha_3\rangle=\pi\sqrt{-1};\quad \langle H,\alpha_1\rangle=\langle H_1,\alpha_2\rangle=\langle
H_1,-\phi\rangle=0.$$ Then $\phi_1=-(\alpha_2+2\alpha_3+2\alpha_4)$
be the maximal root of $\mathfrak k$. Then the Dynkin diagram of
${\mathfrak h}$ is
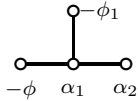
\begin{center}
\setlength{\unitlength}{0.7mm}
\begin{picture}(30,20)(0,-5) \thicklines
\multiput(0,0)(10,0){3}{\circle{2}} \put(10,10){\circle{2}}
\multiput(1,0)(20,0){1}{\line(1,0){8}}
\multiput(11,0)(2,0){1}{\line(1,0){8}} \put(10,1){\line(0,1){8}}
\put(0,-5){\makebox(0,0){$\scriptstyle -\phi$}}
\put(10,-5){\makebox(0,0){$\scriptstyle \alpha_1$}}
\put(15,10){\makebox(0,0){$\scriptstyle -\phi_1$}}
\put(20,-5){\makebox(0,0){$\scriptstyle \alpha_2$}}
\end{picture}
\end{center}
Let ${\mathfrak g}={\mathfrak h}\oplus{\mathfrak
p_1}\oplus{\mathfrak p_2}\oplus{\mathfrak p_3}$ be the decomposition of $\mathfrak g$ corresponding to $(\theta,\tau)$.
Here $\mathfrak k={\mathfrak h}\oplus{\mathfrak p_1}$, which is
a decomposition of $\mathfrak k$ corresponding to the involution
$\tau^{\mathfrak k}$. By Theorem~\ref{subalgebra}, ${\mathfrak p}_1$
is the irreducible representation of ${\mathfrak h}$ with the
highest weight $-\alpha_3$, which is a fundamental dominant weight corresponding
to $\alpha_2$. Let ${\mathfrak k}^{1}=\{x\in {\mathfrak
g}|\theta\tau(x)=x\}={\mathfrak h}\oplus{\mathfrak p}_3$, which is
a decomposition of ${\mathfrak k}^1$ corresponding to the involution
$\tau|_{{\mathfrak k}^1}$. The Dynkin diagram of $F_4$
corresponding to the fundamental system $\{\alpha_1,\alpha_2,\alpha_3'=\alpha_3+\alpha_4,\alpha_4'=-\alpha_4\}$ is
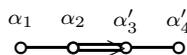
\begin{center}
\setlength{\unitlength}{0.7mm}
\begin{picture}(30,15)(0,-5)
\thicklines \multiput(0,0)(10,0){4}{\circle{2}}
\put(10,0){\circle{2}} \multiput(1,0)(10,0){1}{\line(1,0){8}}
\multiput(11,-0.5)(0,1){2}{\line(1,0){8}} \put(21,0){\line(1,0){8}}
\put(16,-1.5){$>$} \put(0,5){\makebox(0,0){$\scriptstyle \alpha_1$}}
\put(10,5){\makebox(0,0){$\scriptstyle \alpha_2$}}
\put(20,5){\makebox(0,0){$\scriptstyle \alpha_3'$}}
\put(30,5){\makebox(0,0){$\scriptstyle \alpha_4'$}}
\end{picture}
\end{center}
The involution $\theta\tau=e^{\mathrm{ad}(H+H_1)}$ satisfies
$$\langle H+H_1,\alpha_4'\rangle=\pi\sqrt{-1};\quad \langle H+H_1,\alpha_1\rangle=\langle H+H_1,\alpha_2\rangle=\langle H+H_1,\alpha_3'\rangle=0.$$
Here the maximal root $2\alpha_1+3\alpha_2+4\alpha_3'+2\alpha_4'=2\alpha_1+3\alpha_2+4\alpha_3+2\alpha_4=\phi$. It follows that
the Dynkin diagram of ${\mathfrak k}^1$ is
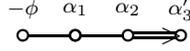
\begin{center}
\setlength{\unitlength}{0.7mm}
\begin{picture}(30,15)(0,-5)
\thicklines \multiput(0,0)(10,0){4}{\circle{2}}
\put(10,0){\circle{2}} \multiput(1,0)(10,0){2}{\line(1,0){8}}
\multiput(21,-0.5)(0,1){2}{\line(1,0){8}} \put(26,-1.5){$>$}
\put(0,5){\makebox(0,0){$\scriptstyle -\phi$}}
\put(10,5){\makebox(0,0){$\scriptstyle \alpha_1$}}
\put(20,5){\makebox(0,0){$\scriptstyle \alpha_2$}}
\put(30,5){\makebox(0,0){$\scriptstyle \alpha_3'$}}
\end{picture}
\end{center}
The involution $\tau=e^{\mathrm{ad}H_1}$ restricted on ${\mathfrak
k}^1$ satisfies
$$\langle H_1,\alpha_3'\rangle=\pi\sqrt{-1};\quad \langle H,\alpha_1\rangle=\langle H_1,\alpha_2\rangle=\langle
H_1,-\phi\rangle=0.$$ Then $-(\alpha_2+2\alpha_3)$ is the maximal
root of $\mathfrak k$, and the Dynkin diagram of ${\mathfrak h}$ is
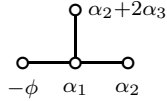
\begin{center}
\setlength{\unitlength}{0.7mm}
\begin{picture}(30,20)(0,-5) \thicklines
\multiput(0,0)(10,0){3}{\circle{2}} \put(10,10){\circle{2}}
\multiput(1,0)(20,0){1}{\line(1,0){8}}
\multiput(11,0)(2,0){1}{\line(1,0){8}} \put(10,1){\line(0,1){8}}
\put(0,-5){\makebox(0,0){$\scriptstyle -\phi$}}
\put(10,-5){\makebox(0,0){$\scriptstyle \alpha_1$}}
\put(20,10){\makebox(0,0){$\scriptstyle \alpha_2+2\alpha_3$}}
\put(20,-5){\makebox(0,0){$\scriptstyle \alpha_2$}}
\end{picture}
\end{center}
By Theorem~\ref{subalgebra}, ${\mathfrak p}_3$ is the irreducible
representation of ${\mathfrak h}$ with the highest weight
$-\alpha_3'$. For the fundamental system $\{-\phi,\alpha_1,\alpha_2,-\phi_1\}$, ${\mathfrak p}_3$ is the irreducible representation
of ${\mathfrak h}_1$ with the highest weight $-(\alpha_1+2\alpha_2+3\alpha_3+\alpha_4)$, which is a
fundamental dominant weight corresponding to $-\phi$. The discussion for $\mathfrak p_2$ is similar. Finally, for the fundamental system $\{-\phi,\alpha_1,\alpha_2,-\phi_1\}$ of ${\mathfrak h}$, we conclude that the highest weights of ${\mathfrak p}_1$, ${\mathfrak p}_2$ and ${\mathfrak p}_3$ as $\mathrm{ad} {\mathfrak h}$ modules are fundamental dominant weights of ${\mathfrak h}$ corresponding to $\alpha_2$, $-\phi_1$ and $-\phi$ respectively, which are pairwise nonisomorphic.
\end{remark}

By the above theory, we classify three-locally-symmetric spaces as follows.
\begin{theorem}\label{tlss}
The classification of three-locally-symmetric spaces $G/H$ for a connected, compact and simple Lie group $G$ is given in Table 1. In Table 1, $A_1=B_1=C_1$, $B_2=C_2$, $A_3=D_3$, $D_1=T$ and $A_0=B_0=C_0=D_0=e$.
\begin{table}[htb]
Table 1: Classification of three-locally-symmetric spaces
\begin{tabular}{lccc|lccc}
\hline
Type & $G$ & $(\theta,\tau^{\mathfrak k})$ & $H$ & Type & $G$ & $(\theta,\tau^{\mathfrak k})$ & $H$  \\
\hline
$A$-I & $A_1$ & (InP1) & $e$ & $A$-II & $A_l$ & (InP2) & $T\times A_{\frac{l-1}{2}}$ \\
&&&&&&& $l\geq 3$ is odd \\
\hline
$A$-III & $A_l$ & (InP3) & $ T^2\times A_{i-1}\times A_{j-i-1}\times A_{l-j}$ &$B$-I & $B_l$ & (InP4) & $B_{l-i}\times D_j\times D_{i-j}$ \\
 &&&  $1\leq i\leq [\frac{l+1}{3}]$ &&&& $2< i\leq l$ \\
 &&& $2i\leq j\leq [\frac{l+i+1}{2}]$ && & & $\frac{i}{2}\leq j\leq i-1$ \\
 \hline
$B$-II & $B_l$ & (InP5) & $B_{i-1}\times B_{l-i}$ & $B$-III & $B_l$ & (InP6) & $B_{i-1}\times B_{i-j}\times B_{l-i}$\\
&&& $\frac{l+1}{2}\leq i\leq l$ & && & $[\frac{2l+3}{3}]\leq i\leq l$\\
&&& &&&& $[\frac{i+2}{2}]\leq j\leq 2i-l$  \\
\hline
$C$-I & $C_l$ &(InP7) & $C_{i}\times C_{j-i}\times C_{l-j}$ & $D$-I & $D_l$ & (InP8) & $D_i\times D_{j-i}\times D_{l-j}$ \\
&&& $1\leq i\leq [\frac{l}{3}]$ &&&& $1\leq i\leq [\frac{l}{3}]$\\
&&& $2i\leq j\leq [\frac{l+i}{2}]$ &&&& $2i\leq j\leq [\frac{l+i}{2}]$\\
\hline
$D$-II & $D_l$ &(InP9) & $B_{i-1}\times D_{l-i}$ & $D$-III & $D_l$ & (InP10) & $D_i\times B_{j-i}\times B_{l-j-1}$ \\
&&& $1\leq i\leq l-2$ &&&& $1\leq i\leq l-2$\\
&&& &&&& $i\leq j\leq [\frac{l+i-1}{2}]$\\
\hline
$D$-IV & $D_l$ & (InP11) & $D_{l-1}$ & $D$-V & $D_l$ & (InP12) &$T^2\times A_{l-2}$  \\
\hline
$E_6$-I & $E_6$ & (InP13) & $T^2\times D_4$ & $E_6$-II & $E_6$ &(InP14) & $T\times A_1\times A_1\times A_3$ \\
\hline
$E_6$-III & $E_6$ & (InP15) & $A_1\times C_{3}$ & $E_7$-I & $E_7$ & (InP16) & $A_1\times A_1\times A_1\times D_4$  \\
\hline
$E_7$-II & $E_7$ &(InP17) & $T\times A_1\times A_5$ & $E_7$-III &$E_7$ & (InP18) & $D_{4}$  \\
\hline
$E_8$-I & $E_8$ & (InP19) & $A_1\times A_1\times D_{6}$ & $E_8$-II & $E_8$ &(InP20) & $D_4\times D_4$  \\
\hline
 $F_4$-I & $F_4$ & (InP21) & $D_4$ & $F_4$-II & $F_4$ &(InP22) & $A_1\times A_1\times C_2$\\
\hline
\end{tabular}
\end{table}
\end{theorem}

\begin{remark}
The well-known examples of three-locally-symmetric spaces are the following:
\begin{enumerate}
\item $SU(2)=SU(2)/\{e\}$,
\item $SU(n_1+n_2+n_3)/S(U(n_1)\times U(n_2)\times U(n_3))$,
\item $SO(l+m+n)/SO(l)\times SO(m)\times SO(n)$,
\item $SO(2n)/U(1)\times U(n-1)$,
\item $Sp(l+m+n)/Sp(l)\times Sp(m)\times Sp(m)$,
\item $E_6/U(1)\times U(1)\times Spin(8)$,
\item $F_4/Spin(8)$.
\end{enumerate}
The first one is the three-locally-symmetric space of type $A$-I in Table 1 of Theorem~\ref{tlss}, the second one is of type $A$-III, the third one corresponds to types $B$-I, $B$-II, $B$-III, $D$-I, $D$-II, $D$-III, and $D$-IV, the fourth one is of type $D$-V, the fifth one is of type $C$-I, the sixth one is of type $E_6$-I, and the seventh one is of type $F_4$-I. For the above cases, $\mathfrak p_1$, $\mathfrak p_2$ and $\mathfrak p_3$ have been proved to be irreducible, and pairwise nonisomorphic with respect to the adjoint action of the Lie algebra $\mathfrak h$ on $\mathfrak p$ except $SO(n+2)/SO(n)$, which is of type $B$-II for $i=l$ and $D$-IV.
\end{remark}

By Theorem~\ref{tlss}, we obtain in Table 2 the dimensions of $\mathfrak p_1$, $\mathfrak p_2$ and $\mathfrak p_3$ for the following cases.
\begin{table}[htb]
Table 2: The dimensions of $\mathfrak p_1$, $\mathfrak p_2$ and $\mathfrak p_3$
\begin{tabular}{cccc|cccc}
\hline
Type & $\dim \mathfrak p_1$ & $\dim \mathfrak p_2$ & $\dim \mathfrak p_3$  & Type & $\dim \mathfrak p_1$ & $\dim \mathfrak p_2$ & $\dim \mathfrak p_3$ \\
\hline
$A$-II & $\frac{(l-1)(l+3)}{4}$ & $\frac{(l+1)(l+3)}{4}$ & $\frac{(l-1)(l+1)}{4}$ & $E_6$-II  & 16 & 16 & 24 \\
\hline
$E_6$-III  & 14 & 28 & 12 & $E_7$-I & 32 & 32 & 32 \\
\hline
$E_7$-II & 24 & 30 & 40 &$E_7$-III& 35 & 35 & 35 \\
\hline
$E_8$-I & 48 & 64 & 64 &$E_8$-II& 64 & 64 & 64 \\
\hline
$F_4$-II & 20 & 8 & 8 \\
\hline
\end{tabular}
\end{table}
Clearly $\mathfrak p_1$, $\mathfrak p_2$ and $\mathfrak p_3$ are pairwise nonisomorphic for types $A$-II, $E_6$-III and $E_7$-II. We can prove the same result for the other cases similar to Remark~\ref{exam}.
In summary, we have the following theorem.
\begin{theorem}\label{thnon}
Let $G/H$ be a three-locally-symmetric space in Theorem~\ref{tlss} with the decomposition $\mathfrak p=\mathfrak p_1\oplus\mathfrak p_2\oplus\mathfrak p_3$. Then $\mathfrak p_1,\mathfrak p_2,\mathfrak p_3$ are pairwise nonisomorphic with respect to the adjoint action of the Lie algebra $\mathfrak h$ on $\mathfrak p$ except types $B$-II for $i=l$ and $D$-IV.
\end{theorem}

\section{Einstein metrics on three-locally-symmetric spaces}
There are some study on the geometry of three-locally-symmetric spaces. In particular, a lot of study on invariant Einstein metrics has been done for some three-locally-symmetric spaces independently. For example,

(1) The flag manifolds $SU(3)/T_{max},Sp(3)/Sp(1)\times Sp(1)\times Sp(1),F_4/Spin(8)$ known as Wallach spaces admit invariant Riemannian metrics of positive section curvature (\cite{Wa1}). The invariant Einstein metrics on the first space are classified in \cite{DN1}, on the other two spaces in \cite{Ro1}. In any case, there exactly four invariant Einstein metrics up to proportionality.

(2) The invariant Einstein metrics on the K$\rm{\ddot{a}}$hler C-spaces $SU(n_1+n_2+n_3)/S(U(n_1)\times U(n_2)\times U(n_3))$, $SO(2n)/U(1)\times U(n-1)$, $E_6/U(1)\times U(1)\times Spin(8)$ are classified in \cite{Ki1}. Every spaces admits four invariant Einstein metrics up to proportionality. Another approach to $SU(n_1+n_2+n_3)/S(U(n_1)\times U(n_2)\times U(n_3))$ is given in \cite{Ar1}.

(3) The Lie group $SU(2)$ considered as $SU(2)/\{e\}$ admits only one left-invariant Einstein metric which is a metric of constant curvature \cite{Be1}.

(4) It is proved in \cite{Ni00} that every three-locally-symmetric space admits at least on invariant Einstein metric. Furthermore, it is proved in \cite{LNF1} that $Sp(l+m+n)/Sp(l)\times Sp(m)\times Sp(n)$ admits exactly four invariant Einstein metrics up to proportionality and that $SO(l+m+n)/SO(l)\times SO(m)\times SO(n)$ admits one, two, three or four invariant Einstein metrics up to proportionality. In particular, it is demonstrated in \cite{Ke1} that $SO(n+2)/SO(n)$ admits just one Einstein metric up to isometry and homothety for $n\geq 3$, the space $SO(4)/SO(2)$ has two such metrics from the classification theorem for five dimensional homogeneous compact Einstein manifolds \cite{ADF1}.

In summary, the classification of invariant Einstein metrics on three-locally-symmetric spaces are well-known for types $A$-I, $A$-III, $B$-I, $B$-II, $B$-III, $C$-I, $D$-I, $D$-II, $D$-III, $D$-IV, $D$-V, $E_6$-I and $F_4$-I in Theorem~\ref{tlss}.

In the following, assume that $G/H$ is a three-locally-symmetric space in Theorem~\ref{tlss} except the above cases. By Theorem~\ref{thnon}, in the decomposition $\mathfrak p=\mathfrak p_1\oplus\mathfrak p_2\oplus\mathfrak p_3$, $\mathfrak p_1,\mathfrak p_2,\mathfrak p_3$ are pairwise nonisomorphic with respect to the adjoint action of the Lie algebra $\mathfrak h$ on $\mathfrak p$. Then we can give the classification of invariant Einstein metrics on these spaces following the theory in \cite{LNF1,Ni00,WZ1}.

Let $d_i$ denote the dimension of $\mathfrak p_i$, and let $\{e_i^j\}$ be an orthonormal basis in $\mathfrak p_i$ with respect to $\langle\cdot,\cdot\rangle=-B(\cdot,\cdot)$, where $i=1,2,3$ and $1\leq j\leq d_i=\dim \mathfrak p_i$. Define the expression $\left[\begin{array}{c} k \\ ij\end{array}\right]$ by
$$\left[\begin{array}{c} k \\ ij\end{array}\right]=\sum_{\alpha,\beta,\gamma}\langle [e_i^\alpha,e_j^\beta],e_k^\gamma\rangle^2,$$
where $\alpha,\beta,\gamma$ range from 1 to $d_i,d_j,d_k$ respectively. Then $\left[\begin{array}{c} k \\ ij\end{array}\right]$ are symmetric in all three indices and $\left[\begin{array}{c} k \\ ij\end{array}\right]=0$ if two indices concide. Let $c_i$ be the Casimir constant of the adjoint representation of $\mathfrak h$ on $\mathfrak p_i$. If $\{e_0^j\}_{1\leq j\leq \dim\mathfrak h}$ is an orthonormal basis in $\mathfrak h$ with respect to $\langle\cdot,\cdot\rangle$ and $e$ is an arbitrary unit vector in $\mathfrak p_i$, then
$$c_i=\sum_j\langle[e_0^j,e],[e_0^j,e]\rangle.$$
For three-locally-symmetric spaces, by \cite{LNF1,WZ1},
\begin{equation}\label{gamma}2A=\left[\begin{array}{c} k \\ ij\end{array}\right]+\left[\begin{array}{c} j \\ ik\end{array}\right]=d_i(1-2c_i).\end{equation}
Let $\rho$ be an invariant metric on $G/H$. We identity it with the corresponding $\mathrm{ad}{\mathfrak h}$-invariant $(\cdot,\cdot)$ on $\mathfrak p$. Since $\mathfrak p_i$ are irreducible and pairwise nonisomorphic, we have
$$(\cdot,\cdot)=x_1\langle\cdot,\cdot\rangle|_{\mathfrak p_1}\oplus x_2\langle\cdot,\cdot\rangle|_{\mathfrak p_2}\oplus x_3\langle\cdot,\cdot\rangle|_{\mathfrak p_3}$$ for some positive real numbers $x_i$. The Ricci curvature $\mathrm{Ric}(\cdot,\cdot)$ of the metric $(\cdot,\cdot)$ is also $\mathrm{ad}\mathfrak h$-invariant. It is easy to see
$$\mathrm{Ric}(\cdot,\cdot)|_{\mathfrak p_i}=r_i(\cdot,\cdot)|_{\mathfrak p_i}$$ for some real numbers $r_i$. As that given in \cite{Ni00}, we have the following formula $$r_i=\frac{1}{2x_i}+\frac{A}{2d_i}(\frac{x_i}{x_jx_k}
-\frac{x_k}{x_ix_j}-\frac{x_j}{x_ix_k}).$$
Here $\{i,j,k\}=\{1,2,3\}$. Put $a_i=\frac{A}{d_i}$. Then we have
\[\left\{ \begin{aligned}r_1=\frac{1}{2x_1}+\frac{a_1}{2}(\frac{x_1}{x_2x_3}
-\frac{x_2}{x_1x_3}-\frac{x_3}{x_1x_2})\\
r_2=\frac{1}{2x_2}+\frac{a_2}{2}(\frac{x_2}{x_1x_3}
-\frac{x_1}{x_2x_3}-\frac{x_3}{x_1x_2})\\
r_3=\frac{1}{2x_3}+\frac{a_3}{2}(\frac{x_3}{x_1x_2}
-\frac{x_1}{x_2x_3}-\frac{x_2}{x_1x_3})\end{aligned}\right.\]
Now the invariant metric $(\cdot,\cdot)$ is Einstein if and only if $r_1=r_2=r_3$. If $a_i=a_j$ for $i\not=j$, then the equations $r_i=r_j$ and $r_i=r_k$ for $k\not=i,j$ become
\[\left\{ \begin{aligned}
&(x_j-x_i)(x_k-2a_i(x_i+x_j))=0,\\
&x_j(x_k-x_i)+(a_i+a_k)(x_i^2-x_k^2)+(a_k-a_i)x_j^2=0
\end{aligned}\right.\]
If $x_j=x_i$, then the second equation is
\begin{equation}\label{equa1}
(1-2a_k)x_i^2-x_ix_k+(a_i+a_k)x_k^2=0.\end{equation}
If $a_k=1/2$, we have only one family of proportional Einstein metrics. Otherwise, $1-2a_k>0$, hence, all real solutions of the equation~(\ref{equa1}) are positive. Then there exist one family of proportional Einstein metrics for $\Delta_1=0$, two families for $\Delta_1>0$, and none for $\Delta_1<0$. Here $\Delta_1$ is the discriminant of (\ref{equa1}). If $x_k=2a_i(x_i+x_j)$, then the second equation is
\begin{equation}\label{equa2}
(a_i+a_k)(1-4a_i^2)x_i^2-(1-2a_i+8a_i^2(a_i+a_k))x_ix_j+(a_i+a_k)(1-4a_i^2)x_j^2=0.
\end{equation}
If $a_i=1/2$ then the equation~(\ref{equa2}) has no solution. Otherwise, $1-4a_i^2>0$, hence, all real roots of the equation~(\ref{equa2}) are positive. Then there exist one family of proportional Einstein metrics for the discriminant $\Delta_2=0$, two families for $\Delta_2>0$, and none for $\Delta_2<0$. Here $\Delta_2$ is the discriminant of (\ref{equa2}).

In particular, if $a_1=a_2=a_3$, then we have the following theorem.
\begin{theorem}[\cite{LNF1} Theorem 3]\label{equiv3}
If $G/H$ is a three-locally-symmetric space in Theorem~\ref{tlss} satisfying $a_1=a_2=a_3$, then, for $a_1\not\in\{\frac{1}{2},\frac{1}{4}\}$, $G/H$ admits exactly four nonproportional invariant Einstein metrics. The parameters $\{x_1,x_2,x_3\}$ has the form $(t,t,t)$, $((1-2a_1)t,2a_1t,2a_1t)$, $(2a_1t,(1-2a_1)t,2a_1t)$, or $(2a_1t,2a_1t,(1-2a_1)t)$. For $a_1=\frac{1}{2}$ and $a_1=\frac{1}{4}$, every invariant Einstein metric is proportional to the standard metric.
\end{theorem}

The following is the method for calculating $c_i$ given in \cite{LNF1}. In detail, $\mathfrak k_i=\mathfrak h\oplus\mathfrak p_i$ is a subalgebra of $\mathfrak g$. Let $K_i$ be the connect Lie subgroup in  $G$ with the Lie algebra $\mathfrak k_i$. In this case, the homogeneous spaces $K_i/H$ and $G/K_i$ are locally symmetric \cite{Be1}. If $K_i$ does not act almost effectively on $M=K/H$, consider its subgroup acting on $M=K_i/H=\widetilde{K_i}/\widetilde{H}$ almost effectively, here $\widetilde{H}$ denotes the corresponding isotropy group. The pair of algebras $(\widetilde{\mathfrak k_i},\widetilde{\mathfrak h})$ is irreducible symmetric \cite{Be1}. If $\widetilde{\mathfrak k_i}$ is simple, then its Killing form $B_{\widetilde{k_i}}$ is proportional to the restriction of the Killing form of $\mathfrak g$ to $\widetilde{\mathfrak k_i}$, i.e., $B_{\widetilde{k_i}}=\gamma_iB|_{\widetilde{k_i}}$. By Lemma 1 in \cite{LNF1}, $$c_i=\gamma_i/2.$$ It follows that $$a_i=\frac{A}{d_i}=\frac{1-\gamma_i}{2}.$$

From the above formulae and results in \cite{LNF1,Ni00,WZ1}, we can classify invariant Einstein metrics on three-locally-symmetric spaces case by case.

\subsection{Invariant Einstein metrics on the three-locally-symmetric space of type $A$-II} For this case, $\mathfrak h\oplus\mathfrak p_2=C_k$, where $l=2k-1$ for $k\geq 2$. By the table for $\gamma_i\geq \frac{1}{2}$ given in \cite{DZ1}, $$\gamma_2=\frac{k+1}{2k}.$$ In fact, there is a method to compute every $\gamma_i$ in \cite{DZ1}. Here for three-locally-symmetric spaces, from (\ref{gamma}) and the dimensions in Table 2, we calculate directly $$\gamma_1=\frac{1}{2},\quad \gamma_3=\frac{k-1}{2k}.$$
It follows that
$$a_1=\frac{1}{4},\quad a_2=\frac{k-1}{4k},\quad a_3=\frac{k+1}{4k}.$$
Let $x_1=1$. The equations $r_1=r_2=r_3$ are equivalent to
\[\left\{ \begin{aligned}
&x_2^2-(2k+1)x_3^2+4kx_2x_3-4kx_2+2k+1=0,\\
&x_2^2-x_3^2+2x_3-2x_2+\frac{1}{k}=0.
\end{aligned}\right.\]
Excluding the summand containing $x_2^2$ form the first equation, we obtain
\begin{equation}\label{A2}x_2(4kx_3-4k+2)=2kx_3^2+2x_3+\frac{1}{k}-2k-1.\end{equation}
For this case, $4kx_3-4k+2\not=0$. Expressing $x_2$ by $x_3$ from $(\ref{A2})$ and inserting it into second one, we obtain
\[\begin{aligned}\label{A22}
&12k^4x_3^4-(48k^4-8k^3)x_3^3+(72k^4-36k^3-4k^2)x_3^2\\
 &-(48k^4-48k^3+4k^2+4k)x_3+12k^4-20k^3+7k^2+2k-1=0.
\end{aligned}\]
Denote by $U_0(x_3)$ the left side of the above equation. Similar to \cite{LNF1}, define Sturm's series by
\[\begin{aligned}
&U_1(x_3)=48k^4x_3^3-(144k^4-24k^3)x_3^2+(144k^4-72k^3-8k^2)x_3-48k^4+48k^3-4k^2-4k,\\
&U_2(x_3)=(6k^3+3k^2)x_3^2-(12k^3-2k^2-\frac{8}{3}k)x_3+6k^3-4k^2-\frac{7}{6}k+\frac{5}{6},\\
&U_3(x_3)=\frac{32k}{27(2k+1)^2}[(108k^4-18k^3-12k^2+4k)x_3-108k^4+36k^3+27k^2-7k-1],\\
&U_4(x_3)=-\frac{27(4k^2-1)^2(27k^4-9k^2+1)}{4(-54k^3+9k^2+6k-2)^2}.
\end{aligned}\]
It follows that
\[\begin{aligned}
&U_0(0)>0, U_1(0)<0, U_2(0)>0, U_3(0)<0, U_4(0)<0;\\
&U_0(\infty)>0, U_1(\infty)>0, U_2(\infty)>0, U_3(\infty)>0, U_4(\infty)<0.
\end{aligned}\]
Denote by $Z(0)$ the number of the sign changes in the series $U_0(0)$, $U_1(0)$, $U_2(0)$, $U_3(0)$, $U_4(0)$ (neglecting zeros) and $Z(\infty)$ the number of sign changes in $U_0(\infty)$, $U_1(\infty)$, $U_2(\infty)$, $U_3(\infty)$, $U_4(\infty)$, where $U_i(\infty)$ denotes the leading coefficient of $U_i(x_3)$ which defines the sign of $U_i(x_3)$ as $x_i\rightarrow\infty$. By Sturm's theorem \cite{Wae1}, the number of real roots of (\ref{A22}) on $(0,\infty)$ is equal to $Z(0)-Z(\infty)=2$ since $U_0(0)\not=0$. By the discussion in \cite{LNF1}, the homogeneous manifold of type $A$-II in Theorem~\ref{tlss} admits exactly two invariant Einstein metrics up to proportionality.
\subsection{Invariant Einstein metrics on the three-locally-symmetric space of type $E_6$-II} For this case, $\mathfrak h\oplus \mathfrak p_1=A_1\oplus A_5$ and $\mathfrak p_1\subset A_5$. By the tables in \cite{DZ1}, $\gamma_1=\frac{1}{2}$. From (\ref{gamma}) and the dimensions in Table 2, $\gamma_2=\frac{1}{2}$ and $\gamma_3=\frac{2}{3}$. It follows that $$a_1=a_2=\frac{1}{4},\quad a_3=\frac{1}{6}.$$ If $x_1=x_2$, then the equation~(\ref{equa1}) is
$$\frac{2}{3}x_1^2-x_1x_3+\frac{5}{12}x_3^2=0.$$
The discriminant $1-\frac{40}{36}<0$, which implies that the equation has no solution. If $x_3=\frac{1}{2}(x_1+x_2)$, then the equation~(\ref{equa1}) equals with
$$15x_1^2-34x_1x_2+15x_2^2=0.$$
The discriminant $34^2-30^2=16^2>0$, which implies that the equation has two solutions.

That is, the homogeneous manifold of type $E_6$-II in Theorem~\ref{tlss} admits exactly two invariant Einstein metrics up to proportionality. The parameters $(x_1,x_2,x_3)$ have the form $(\frac{5}{3}t,t,\frac{4}{3}t)$, or $(\frac{3}{5}t,t,\frac{4}{5}t)$, where $t>0$.

\subsection{Invariant Einstein metrics on the three-locally-symmetric space of type $E_6$-III} For this case, $\mathfrak h\oplus \mathfrak p_1=A_1\oplus A_5$ and $\mathfrak p_1\subset A_5$. By the tables in \cite{DZ1}, $\gamma_1=\frac{1}{2}$. From (\ref{gamma}) and the dimensions in Table 2, $\gamma_2=\frac{3}{4}$ and $\gamma_3=\frac{5}{12}.$ It follows that
$$a_1=\frac{1}{4},\quad a_2=\frac{1}{8},\quad a_3=\frac{7}{24}.$$
Let $x_1=1$. The equations $r_1=r_2=r_3$ are equivalent to
\[\left\{ \begin{aligned}
&x_2^2-13x_3^2+24x_2x_3-24x_2+13=0,\\
&5x_2^2-5x_3^2+12x_3-12x_2+2=0.
\end{aligned}\right.\]
Excluding the summand containing $x_2^2$ form the first equation, we obtain
\begin{equation}\label{E6}x_2(40x_3-36)=20x_3^2+4x_3-21.\end{equation}
For this case, $40x_3-36\not=0$. Expressing $x_2$ by $x_3$ from $(\ref{E6})$ and inserting it into second one, we obtain
$$1200x_3^4-4960x_3^3+7048x_3^2-4152x_3+855=0,$$
which has two real solutions $x_3\approx 1.8845$ or $x_3\approx 0.4838$. From $(\ref{E6})$, $x_2\approx 1.4618$ or $x_2\approx 0.8640$.

That is, the homogeneous manifold of type $E_6$-III in Theorem~\ref{tlss} admits exactly two invariant Einstein metrics up to proportionality. The parameters $(x_1,x_2,x_3)\approx (t,1.4618t,1.8845t)$ or $(x_1,x_2,x_3)\approx (t,0.8640t,0.4838t)$, where $t>0$.

\subsection{Invariant Einstein metrics on the three-locally-symmetric space of type $E_7$-I} For this case, $\mathfrak h\oplus \mathfrak p_1=A_1\oplus D_6$ and $\mathfrak p_1\subset D_6$. By the tables in \cite{DZ1}, $\gamma_1=\frac{5}{9}$. From (\ref{gamma}) and the dimensions in Table 2, $\gamma_2=\gamma_3=\frac{5}{9}.$ It follows that
$$a_1=a_2=a_3=\frac{2}{9}.$$
By Theorem~\ref{equiv3}, i.e. Theorem 3 in \cite{LNF1}, the homogeneous manifold of type $E_7$-I in Theorem~\ref{tlss} admits exactly four invariant Einstein metrics up to proportionality. The parameters $(x_1,x_2,x_3)$ have the form $(t,t,t)$, $(\frac{5}{9}t,\frac{4}{9}t,\frac{4}{9}t)$, $(\frac{4}{9}t,\frac{5}{9}t,\frac{4}{9}t)$, or $(\frac{4}{9}t,\frac{4}{9},\frac{5}{9}t)$, where $t>0$.

\subsection{Invariant Einstein metrics on the three-locally-symmetric space of type $E_7$-II} For this case, $\mathfrak h\oplus \mathfrak p_2=A_1\oplus D_6$ and $\mathfrak p_2\subset D_6$. By the tables in \cite{DZ1}, $\gamma_2=\frac{5}{9}$. From (\ref{gamma}) and the dimensions in Table 2, $\gamma_1=\frac{4}{9}$ and $\gamma_3=\frac{2}{3}$. It follows that
$$a_1=\frac{5}{18},\quad a_2=\frac{2}{9},\quad a_3=\frac{1}{6}.$$
Let $x_1=1$. The equations $r_1=r_2=r_3$ are equivalent to
\[\left\{ \begin{aligned}
&x_2^2+4x_3^2-9x_2x_3+9x_2-4=0,\\
&7x_2^2-7x_3^2+18x_3-18x_2-1=0.
\end{aligned}\right.\]
Excluding the summand containing $x_2^2$ form the first equation, we obtain
\begin{equation}\label{E7}x_2(63x_3-81)=35x_3^2-18x_3-27.\end{equation}
For this case, $63x_3-81\not=0$. Expressing $x_2$  by $x_3$ from $(\ref{E7})$ and inserting it into second one, we obtain
$$2744x_3^4-13482x_3^3+24732x_3^2-19926x_3+5832=0,$$
which has two real solutions $x_3\approx 1.5535$ or $x_3\approx 0.7302$. From $(\ref{E7})$, $x_2\approx 1.7489$ or $x_2\approx 0.6139$.

That is, the homogeneous manifold of type $E_7$-II in Theorem~\ref{tlss} admits exactly two invariant Einstein metrics up to proportionality. The parameters $(x_1,x_2,x_3)\approx (t,1.7489t,1.5535t)$ or $(x_1,x_2,x_3)\approx (t,0.6139t,0.7302t)$, where $t>0$.

\subsection{Invariant Einstein metrics on the three-locally-symmetric space of type $E_7$-III} For this case, $\mathfrak h\oplus \mathfrak p_1=A_7$. It is the same as that of type $E_7$-II. That is, $\gamma_1=\frac{4}{9}$. From (\ref{gamma}) and the dimensions in Table 2, $\gamma_2=\gamma_3=\frac{4}{9}$. It follows that
$$a_1=a_2=a_3=\frac{5}{18}.$$
By Theorem~\ref{equiv3}, i.e. Theorem 3 in \cite{LNF1}, the homogeneous manifold of type $E_7$-III in Theorem~\ref{tlss} admits exactly four invariant Einstein metrics up to proportionality. The parameters $(x_1,x_2,x_3)$ have the form $(t,t,t)$, $(\frac{4}{9}t,\frac{5}{9}t,\frac{5}{9}t)$, $(\frac{5}{9}t,\frac{4}{9}t,\frac{5}{9}t)$, or $(\frac{5}{9}t,\frac{5}{9}t,\frac{4}{9}t)$, where $t>0$.

\subsection{Invariant Einstein metrics on the three-locally-symmetric space of type $E_8$-I} For this case, $\mathfrak h\oplus \mathfrak p_2=A_1\oplus E_7$ and $\mathfrak p_2\subset E_7$. By the tables in \cite{DZ1}, $\gamma_2=\frac{3}{5}$. From (\ref{gamma}) and the dimensions in Table 2, $\gamma_1=\frac{7}{15}$ and $\gamma_3=\frac{3}{5}.$
It follows that
$$a_1=\frac{4}{15},\quad a_2=a_3=\frac{1}{5}.$$
If $x_2=x_3$, then the equation~(\ref{equa1}) is
$$\frac{7}{15}x_2^2-x_1x_2+\frac{7}{15}x_1^2=0.$$
The discriminant $1-(\frac{14}{15})^2>0$, which implies that the equation has two solutions. If $x_1=\frac{2}{5}(x_2+x_3)$, then the equation~(\ref{equa1}) equals with
$$147x_2^2-281x_2x_3+147x_3^2=0.$$
The discriminant $281^2-294^2<0$, which implies that the equation has no solution.

That is, the homogeneous manifold of type $E_8$-I in Theorem~\ref{tlss} admits exactly two invariant Einstein metrics up to proportionality. The parameters $(x_1,x_2,x_3)$ have the form $(qt,t,t)$, where $t>0$ and $q$ is the root of the equation $7x^2-15x+7=0$.

\subsection{Invariant Einstein metrics on the three-locally-symmetric space of type $E_8$-II} For this case, $\mathfrak h\oplus \mathfrak p_1=D_8$. It is the same as that of type $E_8$-I. That is, $\gamma_1=\frac{7}{15}$. From (\ref{gamma}) and the dimensions in Table 2, $\gamma_2=\gamma_3=\frac{7}{15}.$
It follows that
$$a_1=a_2=a_3=\frac{4}{15}.$$
By Theorem~\ref{equiv3}, i.e. Theorem 3 in \cite{LNF1}, the homogeneous manifold of type $E_8$-II in Theorem~\ref{tlss} admits exactly four invariant Einstein metrics up to proportionality. The parameters $(x_1,x_2,x_3)$ have the form $(t,t,t)$, $(\frac{7}{15}t,\frac{8}{15}t,\frac{8}{15}t)$, $(\frac{8}{15}t,\frac{7}{15}t,\frac{8}{15}t)$, or $(\frac{8}{15}t,\frac{8}{15}t,\frac{7}{15}t)$, where $t>0$.

\subsection{Invariant Einstein metrics on the three-locally-symmetric space of type $F_4$-II} For this case, $\mathfrak h\oplus \mathfrak p_1=D_4$. By the tables in \cite{DZ1}, $\gamma_1=\frac{7}{9}$. From (\ref{gamma}) and the dimensions in Table 2, $\gamma_2=\gamma_3=\frac{4}{9}.$
It follows that
$$a_1=\frac{1}{9},\quad a_2=a_3=\frac{5}{18}.$$
If $x_2=x_3$, then the equation~(\ref{equa1}) is
$$\frac{7}{9}x_2^2-x_1x_2+\frac{7}{18}x_1^2=0.$$
The discriminant $1-\frac{98}{81}<0$, which implies that the equation has no solution.  There exists none Einstein metrics. If $x_1=\frac{5}{9}(x_2+x_3)$, then the equation~(\ref{equa1}) equals with
$$196x_2^2-499x_2x_3+196x_3^2=0.$$
The discriminant $499^2-392^2>0$, which implies that the equation has two solutions.

That is, the homogeneous manifold of type $F_4$-II in Theorem~\ref{tlss} admits exactly two invariant Einstein metrics up to proportionality. The parameters $(x_1,x_2,x_3)$ have the form $(\frac{5}{9}(q+1)t,qt,t)$, where $t>0$ and $q$ is the root of the equation $196x^2-499x+196=0$.

\section{Acknowledgments}
This work is supported by National Natural Science Foundation of
China (No. 11001133). The first author would like to thank Prof.
J.A. Wolf for the helpful conversation and suggestions.

\end{document}